\documentclass[10pt,a4paper]{article}
\usepackage{amsmath}
\usepackage{amsfonts}
\usepackage{amssymb}
\usepackage{a4,a4wide}
\usepackage{graphics}
\usepackage{color}
\usepackage[T1]{fontenc}

\def\ep{{\varepsilon}}
\def\R{\mathbb R}
\def\c{c^*}

\newcommand{\fdem}{ \hfill $\square$ }

\newtheorem{theo}{\textbf{Theorem}}[section]

\newtheorem{lem}[theo]{\textbf{Lemma}}
\newtheorem{prop}[theo]{\textbf{Proposition}}

\newtheorem{rem}[theo]{\textbf{Remark}}

\title{Bistable travelling waves for nonlocal reaction diffusion equations}
\date{}

\begin{document}

\maketitle

\begin{center}
{\large\bf Matthieu Alfaro \footnote{ I3M, Universit\'e de
Montpellier 2, CC051, Place Eug\`ene Bataillon, 34095 Montpellier
Cedex 5, France. E-mail: malfaro@math.univ-montp2.fr},
J{\'e}r{\^o}me Coville \footnote{Equipe BIOSP, INRA Avignon,
Domaine Saint Paul, Site Agroparc, 84914 Avignon Cedex 9, France.
E-mail: jerome.coville@avignon.inra.fr} and  Ga\"el Raoul
\footnote{Centre d'\'Ecologie Fonctionnelle et \'Evolutive, UMR
5175, CNRS, 1919 Route de Mende, 34293 Montpellier, France.
E-mail: raoul@cefe.cnrs.fr
.}.}\\
[2ex]

\end{center}



\vspace{10pt}

\begin{abstract}
 We are concerned with travelling wave solutions
arising in a reaction diffusion equation with bistable and
nonlocal nonlinearity, for which the comparison principle does not
hold. Stability of the equilibrium $u\equiv 1$ is not assumed. We
construct a travelling wave solution connecting 0 to an unknown
steady state, which is \lq\lq above and away\rq\rq\, from the
intermediate equilibrium. For focusing kernels we prove that, as
expected, the wave connects 0 to 1. Our results also apply readily
to the nonlocal ignition case.
\\


\end{abstract}

\section{Introduction}\label{s:intro}

We consider the nonlocal bistable reaction diffusion equation
\begin{equation}\label{edp}
\partial _t u=\partial _{xx} u +u(u-\theta)(1-\phi
* u)\quad \text{ in }(0,\infty)\times \R,
\end{equation}
where $0<\theta<1$. Here $\phi*u(x):=\displaystyle \int _\R
u(x-y)\phi (y)\,dy$, with $\phi$  a given bounded kernel such that
\begin{equation}\label{noyau}
\phi \geq 0,\quad \phi(0)>0,\quad \int _\R \phi =1.
\end{equation}
We are looking for travelling waves solutions supported by the
integro-differential equation (\ref{edp}), that is a speed $\c \in
\R$ and a smooth $U$ such that
\begin{equation}\label{pb-tw}
-U''-\c U'=U(U-\theta)(1-\phi * U)\quad \text{ in }\R,
\end{equation}
supplemented with the   boundary conditions
\begin{equation}\label{weak}
\liminf_{x\to -\infty} U(x)>\theta,\quad \lim_{x\to+\infty}U(x)=0.
\end{equation}
In this work we construct such a travelling wave solution, and
show that the behavior on the left is improved to
$\lim_{x\to-\infty}U(x)=1$ for {\it focusing kernels}. Our results
also apply readily to the nonlocal {\it ignition case} (see
below).

\medskip

Equation \eqref{edp} is a nonlocal version of the well known
reaction-diffusion equation
\begin{equation}\label{rd}
\partial _t u=\partial _{xx} u +f(u)\quad \text{ in }(0,\infty)\times \R,
\end{equation}
with the bistable nonlinearity $f(s):=s(s-\theta)(1-s)$.
Homogeneous reaction diffusion equations have been extensively
studied in the literature (see \cite{Kan1}, \cite{Aro-Wei1,
Aro-Wei2}, \cite{Fif-Mac}, \cite{Ber-Nic-Sch},
\cite{Volpert-Volpert-Volpert} among others) and are  known to
support the existence of monotone travelling fronts for three
classes of nonlinearity: bistable, ignition  and monostable.
Moreover, for bistable and ignition nonlinearities, there exists a
unique front speed $c^*$ whereas, for monostable nonlinearities,
there exists a critical speed $c^*$ such that all speeds $c\geq
c^*$ are admissible. In this local context, many techniques based
on the comparison principle --- such as some monotone iterative
schemes or the sliding method \cite{Ber-Nir3}--- can be used to
get {\it a priori} bounds, existence and monotonicity properties
of  the fronts.

Recently, much attention was devoted to the introduction of a
nonlocal effect into the nonlinear reaction term. From the
mathematical point of view, the analysis is quite involved since
integro-differential equations with a nonlocal competition term
generally do not satisfy the comparison
 principle.  In \cite{Ber-Nad-Per-Ryz}, Berestycki, Nadin, Perthame and Ryzhik
have  considered the following non-local version of the Fisher-KPP
equation
\begin{equation}\label{kppnonloc}
\partial _t u=\partial _{xx} u +u(1-\phi* u)\quad \text{ in }(0,\infty)\times \R.
\end{equation}
They  prove that equation \eqref{kppnonloc} admits a critical
speed $c^*$ so that,  for all $c \ge c^*$, there exists a
travelling wave $(c,U)$ solution of
\begin{equation*}
\begin{cases}
-U''-cU'=U(1-\phi* U)\quad \text{ in } \R.\\
\liminf _{x\to -\infty} U(x)>0, \quad \lim_{x\to+\infty}U(x)=0.
\end{cases}
\end{equation*}
In favorable situations, namely when the steady state $u\equiv 1$
remains linearly stable, they further obtain
$\lim_{x\to-\infty}U(x)=1$. Nevertheless, the positive steady
state $u\equiv 1$ may present, for some kernels, a Turing
instability (see e.g. \cite{Gen-Vol-Aug},  \cite{Ber-Nad-Per-Ryz},
\cite{Alf-Cov}). In such situations, it was proved in
\cite{Fan-Zha} and in \cite{Alf-Cov} that the waves with large
speeds actually connect the two unstable states 0 and 1. Notice
that the former work considers kernels with exponential decay and
uses monotonicity arguments inspired by \cite{Gom-Tro}, whereas
the latter uses more direct arguments which allow kernels with
algebraic decay. Concerning this issue of the behavior of the wave
on the left, we also refer the reader to \cite{Nad-Per-Tan},
\cite{Nad-Per-Ros-Ryz}. In a related framework, the authors of the
present work have constructed curved fronts for nonlocal reaction
diffusion equations \cite{A-Cov-Rao} of the form
$$
\begin{aligned}
\partial _t u(t,x,y)=&\Delta_{xx} u(t,x,y) +\partial_{yy}u(t,x,y)+\\
&u(t,x,y)\left(r(v-Bx.e)-\int_{\R}k(y-Bx.e,y'-Bx.e)u(t,x,y')\,dy'\right)
\end{aligned}
$$
for $t>0$, $x\in \R ^d$ (spatial variable), $y\in \R$
(phenotypical trait). In population dynamics, such equations serve
as prototypes of models for structured populations evolving in a
environmental cline.

\medskip

In view of the existence of fronts for both the nonlocal
Fisher-KPP equation \eqref{kppnonloc} and the local bistable
equation \eqref{rd}, it is then expected that the nonlocal
bistable equation \eqref{edp} supports the existence of travelling
waves. In this work, we shall construct such a solution.  It is
worth being mentioned that, among other things, nonlinearities
such as $u(\phi*u-\theta)(1-u)$ are treated in \cite{Wan-Li-Ru}.
Notice that our equation does not fall into \cite[equation
 (1.6)]{Wan-Li-Ru} since $g(u,v)=u(1-u)(1-v)$ does not satisfy
 \cite[hypothesis (H1)]{Wan-Li-Ru}, which actually provides the stability
 of both  $u\equiv 0$ and $\equiv 1$.

Let us now state our main result on the existence of a travelling
wave solution.

\begin{theo}[A bistable travelling wave]\label{th:tw} There exist a speed $\c\in \R$ and a positive profile $U\in C^2(\R )$ solution of
\begin{equation}\label{eq-dans-espace}
-U''-\c U'= U(U-\theta)(1-\phi * U) \quad \text{ on } \R ,
\end{equation}
such that, for some $\ep >0$,
\begin{equation}\label{gauche}
U(x)\geq \theta +\ep\; \text{ for all } x\in(-\infty,-1/\ep),
\end{equation}
$U$ is decreasing on $[\bar x,+\infty)$ for some $\bar x>0$, and
\begin{equation}\label{droite}
 \lim_{x\to+\infty}U(x)=0.
\end{equation}
\end{theo}

Now, if the kernel tends to the Dirac mass, we expect the above
travelling wave to be a perturbation of the underlying wave for
the local case, namely $(\c_0,U_0)$ the unique solution of
\begin{equation}\label{local}
\begin{cases}
{U_0}''+\c_0{U_0}'+U_0(U_0-\theta)(1-U_0)=0,\\
\lim_{x\to-\infty}U_0(x)=1, \quad U_0(0)=\theta, \quad
\lim_{x\to+\infty}U_0(x)=0,
\end{cases}
\end{equation}
and so to satisfy $\lim_{x\to-\infty}U(x)=1$. Our next result
states such a behavior assuming $\c _0\neq 0$, which is equivalent
to $\theta \neq \frac 12$. To make this perturbation analysis
precise, we take $\sigma >0$ as a focusing parameter, define
\begin{equation}\label{noyau-focusing}
\phi_\sigma(x):=\frac 1 \sigma \phi \left(\frac x \sigma\right),
\end{equation}
and are interested in the asymptotics $\sigma \to 0$.

\begin{prop}[Focusing kernels]\label{prop:focusing} Denote by $(\c _\sigma, U_\sigma)$
the travelling wave associated with the kernel $\phi _\sigma$, as
constructed in Theorem \ref{th:tw}.
\begin{itemize}
\item [(i)] Assume $\int _\R |z|\phi (z)\,dz<\infty$. Then $\c
_\sigma \to \c _0$, as $\sigma \to 0$. \item [(ii)] Assume $\theta
\neq \frac 12$ and $\int _{\R} z^2 \phi(z)\,dz<\infty$.  Then
there is $\sigma _0>0$ such that, for all $0<\sigma<\sigma _0$,
$$
\lim_{x\to-\infty}U_\sigma(x)=1.
$$
\end{itemize}
\end{prop}

\begin{rem}[Ignition case] While proving the above results for the
bistable case, it will become clear that the same  (with the
additional information $\c>0$) holds for the ignition case, that
is
$$ -U''-\c U'=\begin{cases} 0 &
\text{ where } U<\theta\\
(U-\theta)(1-\phi *U) &\text{ where } U\geq \theta,
\end{cases}
$$
for which proofs are simpler. This will be clarified in Section
\ref{s:ignition}.
\end{rem}

Let us comment on the main result, Theorem \ref{th:tw}. Due to the
lack of comparison principle, the construction of a travelling
wave solution is based on  a topological degree argument, a method
introduced initially in \cite{Ber-Nic-Sch}. After establishing a
series of {\it a priori} estimates, it enables to construct a
solution in a bounded box. Then we let the size of the box tend to
infinity to construct a solution on the whole line. In contrast
with \cite{Ber-Nad-Per-Ryz} and because of the intermediate
equilibrium $u\equiv \theta$, it is far from clear that the
constructed wave is non trivial --- or, equivalently, that it
\lq\lq visits\rq\rq\, both $(0,\theta)$ and $(\theta,1)$. Such an
additional difficulty also arises in the construction of bistable
waves in cylinders \cite{berestycki-nirenberg}, where the authors
use energy arguments to exclude the possibility of triviality.
This seems not to be applicable to our nonlocal case. Our
arguments are rather direct and are based on the sharp property of
Proposition \ref{u0iffx0} and the construction of bump-like
sub-solutions in Lemma \ref{functionphi}.

\medskip

The organization of the paper is as follows. In Section
\ref{s:tw}, we construct a solution $u$ on a bounded interval
thanks to a Leray-Schauder topological degree argument. In Section
\ref{s:nontrivial}, we show that, when we let the bounded interval
tend to the whole line, the limit profile $U$ is non trivial. The
behaviors \eqref{gauche} and \eqref{droite} are then proved in
Section \ref{s:properties}. We  investigate the case of the
focusing kernels, that is Proposition \ref{prop:focusing}, in
Section \ref{s:focusing}. Last, in Section \ref{s:ignition} we
indicate how to handle the ignition case.

\section{Construction of a solution $u$ in a box}\label{s:tw}

Notice that the methods used in this section are inspired by
\cite{Ber-Nad-Per-Ryz}.

\medskip

For $a>0$ and $0\leq \tau \leq 1$, we consider the problem of
finding a speed $c=c_\tau^a\in \R$ and a profile
$u=u_\tau^a:[-a,a]\to \R$ such that
$$
P_\tau(a)\quad\begin{cases}
\,-u''-cu^{\prime}=\tau \mathbf{1} _{\{u\geq 0\}} u(u-\theta)(1-\phi * \bar u)\quad \text{ in }(-a,a)\vspace{5pt}\\
\,u(-a)=1,\qquad u(0)=\theta,\qquad  u(a)=0,\\
\end{cases}
$$
where $\bar u$ denotes the extension of $u$ equal to $1$ on
$(-\infty,-a)$ and 0 on $(a,\infty)$ (in the sequel, for ease of
notation, we always write $u$ in place of $\bar u$). This realizes
a homotopy from a local problem ($\tau =0$) to our nonlocal
problem ($\tau=1$) in the box $(-a,a)$. We shall construct a
solution to $P_1(a)$ by using a Leray-Schauder topological degree
argument.

If $(c,u)$ is a solution achieving a negative minimum at $x_m$
then $x_m\in(-a,a)$ and $-u''-cu'=0$ on a neighborhood of $x_m$.
The maximum principle thus implies $u\equiv u(x_m)$, which cannot
be. Therefore any solution  of $P_\tau(a)$ satisfies $u\geq 0$
and, by the strong maximum principle,
\begin{equation}
u> 0\quad\text{and }\quad  -u''-cu^{\prime}= \tau
u(u-\theta)(1-\phi
*  u)\quad \text{ in }(-a,a). \label{alternative}
\end{equation}

\subsection{A priori bounds of solutions in the box}\label{ss:box}

The following lemma provides {\it a priori} bounds for $u$.

\begin{lem}[A priori bounds for $u$]\label{lem:a-priori-bounds}
 There exist $M>1$ (depending only on the kernel $\phi$) and $a_0>0$ such that, for all $a\geq a_0$ and all $0\leq \tau\leq 1$, any
 solution $(c,u)$ of $P_\tau(a)$ satisfies
 $$
  0\leq u(x) \leq M,\quad \forall x\in[-a,a].
  $$
\end{lem}

\noindent {\bf Proof.}  If $\tau =0$ we directly get $0\leq u\leq
1$ for the local problem. Now, for $0<\tau \leq 1$, assume
$M:=\max_{x\in[-a,a]} u(x)>1$ (otherwise there is nothing to
prove). In view of the boundary conditions, there is a
$x_m\in(-a,a)$ such that $M=u(x_m)$. Evaluating
\eqref{alternative} at $x_m$ we see that $\phi*  u(x_m)\le 1$.

Now since $ u\geq 0$, we also have
\begin{equation}\label{majc.eq.2}
 -u^{\prime\prime}-cu^{\prime}=\tau u(u-\theta)(1-\phi*  u)\leq u^2 + \theta u (\phi * u)\leq (1+\theta)M^2\leq 2 M^2.
\end{equation}

Let us first assume that $c<0$. For $x\in[-a,x_m]$ it follows from
\eqref{majc.eq.2} that
$$
\int_{x}^{x_m}\left(u^{\prime}(z)e^{-|c|z}\right)^{\prime}dz\ge
-\int_{x}^{x_m}2M^2e^{-|c|z}\, dz.$$ Using $u'(x_m)=0$, isolating
$u'(x)$ and integrating again from $x$ to $x_m$, we discover
$$
\int_{x}^{x_m}u^{\prime}(z)\,dz\le
-\frac{2M^2}{|c|}\int_{x}^{x_m}(e^{-|c|(x_m-z)}-1)\,dz.
$$
Using $u(x_m)=M$ and isolating $u(x)$, we get after elementary
computations
$$
u(x)\ge M\left[1-2M(x-x_m)^2 B(|c|(x_m-x))\right],
$$
where $B(y):=\frac{e^{-y}+y-1}{y^2}$. Observe that $B(y)\le
\frac{1}{2}$ for $y>0$ so that
\begin{equation}
u(x)\ge M(1-M(x-x_m)^2),\quad \forall x\in[-a,x_m],
\label{majc.eq.3}
\end{equation}
and in particular, for $x=-a$,
\begin{equation}
1\ge M(1-M(a+x_m)^2). \label{majc.eq.4}
\end{equation}

Now we define
\begin{equation}\label{choix}
x_0:=\frac 1 {2\sqrt M}.
\end{equation}
If $x_m\in(-a,-a+x_0]$, then \eqref{majc.eq.4} shows that
$M\leq\left(1-M{x_0}^2\right)^{-1}=\frac{4}{3}$. If
$x_m\in[-a+x_0,a)$, then
$$
1\ge \phi*  u(x_m)\ge
\int_{0}^{x_0}\phi(z)u(x_m-z)\,dz \ge M
\int_{0}^{x_0}\phi(z)(1-Mz^2)\,dz,
$$
where we have used \eqref{majc.eq.3}. From the definition of $x_0$
we deduce that
$$
1\geq \frac 34 M \int _0 ^{1/2\sqrt M}\phi(z)\,dz\geq \frac 34
M\left(\int _0 ^{1/2\sqrt M}(\phi (0)-\Vert \phi
'\Vert_{L^\infty(-1,1)}z)\,dz\right),
$$
so that
\begin{equation}\label{borne}
M\leq \left(\frac{8}{3}\frac{1+\frac 3{32}\Vert \phi
'\Vert_{L^\infty(-1,1)}}{\phi(0)}\right)^2.
\end{equation}
This concludes the proof in the case $c<0$. The case $c>0$ can be
treated in a similar way by working on $[x_m,a]$ rather than on
$[-a,x_m]$.

Last if $c=0$, by integrating twice the inequality $-u''\leq 2M^2$
on $[x,x_m]$ we directly obtain \eqref{majc.eq.3} and we can
repeat the above arguments. This completes the proof of the lemma.
\fdem

\medskip
We now provide {\it a priori} bounds for the speed $c$.

\begin{lem}[A priori upper bound for $c$]\label{lem:a-priori-speed}
There exists $a_0>0$ such that, for all $a\geq a_0$ and all $0\leq
\tau\leq 1$, any solution $(c,u)$ of $P_\tau(a)$ satisfies $c\leq
2\sqrt{2M}=:c_{max}$, where $M$ is the upper bound for $u$ defined
in Lemma \ref{lem:a-priori-bounds}.
\end{lem}

\noindent {\bf Proof.} Since $ -u^{\prime\prime}-cu^{\prime}\leq
u^2 + \theta u (\phi * u)
 \leq (1+\theta)Mu\leq 2 Mu$, we can reproduce the proof of
 \cite[Lemma 3.2]{Ber-Nad-Per-Ryz} with $\mu \leftarrow 2M$.\fdem

 \medskip

We now provide {\it a priori} bounds for the speed $c$. We will
prove two separate estimates.

\begin{lem}[A priori lower bound for $c$, uniform w.r.t. $\tau$]\label{lem:a-priori-speed3}
For any $a>0$, there exists  $\tilde c_{min}(a) >0$ such that, for
all $0\leq \tau \leq 1$, any solution $(c,u)$ of $P_\tau(a)$
satisfies $c\geq -\tilde c_{min}(a)$.
\end{lem}

\noindent {\bf Proof.} Let $a>0$ be given. We consider a solution
$(c,u)$ of $P_\tau(a)$. It satisfies:
$$-u''-cu'+(M^2+1)u\geq 0,$$
as well as $u(-a)=1$ and $u(a)=0$. Since $M^2+1\geq 0$, the
comparison principle applies and $u\geq v$, where $v$ is the
solution of $-v''-cv'+(M^2+1)v=0$ such that $v(-a)=1$ and
$v(a)=0$. Explicitly computing $v$, we get
$$
v(0)=\frac{1-e^{(\lambda _+ - \lambda _-)a}}{e^{-\lambda _+
a}-e^{(\lambda _+-2\lambda _-)a}}\,,\quad \lambda _ \pm := \frac
{-c\pm \sqrt{c^2+4(M^2+1)}}2.
$$
We see that $v(0)\to 1$ as $c\to -\infty$. It follows that, for
any $a>0$, there exists $\tilde c_{min}(a) >0$ such that $c\leq
-\tilde c_{min}(a)$ implies $\theta<v(0)\leq u(0)$, so that $u$
cannot solve $P_\tau(a)$. Hence, all solutions $(c,u)$ of
$P_\tau(a)$ with $0\leq \tau \leq 1$ are such that $c\geq -\tilde
c_{min}(a)$. \fdem

\medskip


 \begin{lem}[A priori lower bound for $c$ for $\tau =1$, uniform w.r.t. $a$]\label{lem:a-priori-speed2}
There exist  $c_{min} >0$ and $a_0>0$ such that, for all $a\geq
a_0$, any solution $(c,u)$ of $P_1(a)$ satisfies $c\geq -c_{min}$.
\end{lem}

\noindent {\bf Proof.} We assume that $c\leq -1$ (otherwise there
is nothing to prove), and define $M
>0$ and $a_0>0$ as in Lemma \ref{lem:a-priori-bounds}.

The first step of the proof is to find uniform bounds on $u'$,
following ideas from \cite[Lemma 3.3]{Ber-Nad-Per-Ryz}. We first
notice that $(e^{cx}u'(x))'=e^{cx}\left(u''(x)+cu'(x)\right)$, and
an integration of this expression provides for $x>y$:
$$e^{cx}u'(x)-e^{cy}u'(y)=\int_y^x e^{cz} u(z)(u(z)-\theta)(1-\phi\ast u(z))\,dz.$$
Thank to Lemma \ref{lem:a-priori-bounds}, we have $\left|
u(u-\theta)(1-\phi\ast u)\right|\leq M(M+\theta)(1+M)=:Q$, so that
\begin{equation}\label{inegalite-nadin}
u'(y)e^{|c|(x-y)}-\frac {Q}{|c|} e^{|c|(x-y)}\leq u'(x)\leq
u'(y)e^{|c|(x-y)}+\frac {Q}{|c|} e^{|c|(x-y)},\quad \forall x,y\in
[-a,a],\,x>y,
\end{equation}
\begin{equation}\label{derivee-par-dessus}
 u'(y) \leq \frac{2 Q}{|c|},\quad \forall y\in(-a,a),
\end{equation}
where we have chosen $x=a$, and used the fact that $u'(a)\leq 0$
to obtain this last estimate.

Next, define
$$
K_0:=2 \max _{c\leq -1}\frac 1 {|c|} \ln \left(\frac{M c^2}
Q+1\right).
$$
We claim that, for all $c\leq -1$, all $a\geq a_0$,
\begin{equation}\label{derivee-par-dessous}
-\frac{2 Q}{|c|} \leq u'(x),\quad \forall x\in(-a,a-K_0].
\end{equation}
Indeed, assume by contradiction that there are some $c\leq -1$,
$a\geq a_0$, $y\in(-a,a-K_0]$ such that $u'(y)<-\frac{2Q}{|c|}$.
{}From \eqref{inegalite-nadin} we deduce that $u'(x)\leq -\frac
Q{|c|} e^{|c|(x-y)}$ for $x>y$. Integrating this from
 $y$ to $a$ and using $u(a)=0$ we see that
$$
M\geq u(y) \geq  \frac Q{c^2}(e^{|c|(a-y)}-1) \geq \frac
Q{c^2}(e^{|c|K_0}-1),
$$
which contradicts the definition of $K_0$. This proves
\eqref{derivee-par-dessous}.

Next, since $\phi\in L^1(\mathbb R)$, there exists $R>0$ such that
$M\int_{[-R,R]^c}\phi\leq \frac{1-\theta}8$. Thanks to the
conditions $u(-a)=1$ and $u(0)=\theta$ in $P_1(a)$, we can define
$x_0<0$ as the  largest negative real such that
$u(x_0)=\theta+\frac{1-\theta}2$. We can use
\eqref{derivee-par-dessous} to estimate $u(x)$ from below for
$x\in[x_0-R,x_0+2R]\cap [-a,a]$:
\begin{equation}
 u(x)\geq\theta+\frac{1-\theta}2-\frac{2Q}{|c|}2R\geq\theta+\frac{1-\theta}4,\label{borneinf}
\end{equation}
as soon as $c\leq -\frac{16QR}{1-\theta}$. Similarly, using
\eqref{derivee-par-dessus},
\begin{equation}\label{bornesup}
u(y)\leq \theta+\frac{1-\theta}2+\frac{2Q}{|c|}2R\leq
\theta+\frac{3(1-\theta)}4,
\end{equation}
as soon as $c\leq -\frac{16QR}{1-\theta}$. In particular
\eqref{borneinf} and \eqref{bornesup} imply that
$[x_0-R,x_0+2R]\subset (-a,0)$ if $-c$ is large enough.
We then estimate $\phi\ast u(x)$ for $x\in [x_0,x_0+R]$:
\begin{eqnarray}
 \phi\ast u(x)&\leq& \int_{[-R,R]}\phi(y)u(x-y)\,dy+\int_{[-R,R]^c}\phi(y)u(x-y)\,dy\nonumber\\
 &\leq&\max_{[x_0-R,x_0+2R]}u+M\int_{[-R,R]^c}\phi\nonumber\\
 &\leq & \theta+\frac {1-\theta}2+\frac{2Q}{|c|}2R+\frac{1-\theta}8\nonumber\\
 &\leq& 1-\frac{1-\theta}8,\label{estimphi}
\end{eqnarray}
as soon as $c\leq -\frac{16QR}{1-\theta}$.

If $u$ is not non-increasing on $[x_0,x_0+R]$, the definition of
$x_0$ implies the existence of a local minimum $\bar x\in
(x_0,x_0+R)$. An evaluation \eqref{alternative} in $\bar x$ then
shows that $1\leq \phi\ast u(\bar x)$, which is only possible if
$-c$ is not too large, thanks to \eqref{estimphi}.

If on the contrary, $u$ is non-increasing on $[x_0,x_0+R]$ and
$c\leq 0$, then, for $[x_0,x_0+R]$,
$$u''(x)\leq u''(x)+cu'(x)=-u(x)(u(x)-\theta)(1-\phi\ast u(x))\leq -\theta \frac{(1-\theta)^2}{32}.$$
It follows that $u'(x_0)-u'(x_0+R)\geq
\frac{\theta(1-\theta)^2}{32}R$, which, combined to
\eqref{derivee-par-dessous} and \eqref{derivee-par-dessus} implies
$$c\geq -\frac{128Q R}{\theta(1-\theta)^2},$$
so that in any case, $c_{min}:=-\frac{128Q R}{\theta(1-\theta)^2}$
is an explicit lower bound for $c$.\fdem

\subsection{Construction of a solution in the box}\label{ss:construction}

Equipped with the above {\it a priori} estimates, we now use a
Leray-Schauder topological degree argument (see e.g.
\cite{Ber-Nic-Sch}, \cite{Ber-Nad-Per-Ryz} or \cite{A-Cov-Rao} for
related arguments) to construct a solution $(c,u)$ to $P_1(a)$.

\begin{prop}[A solution in the
box]\label{prop:sol-boite} There exist $K>0$ and $a_0>0$ such
that, for all $a\geq a_0$, Problem $P_1(a)$ admits a solution
$(c,u)$, that is
$$
\begin{cases}
\,-u''-cu^{\prime}=u(u-\theta)(1-\phi *  u)\quad \text{ in }(-a,a)\vspace{3pt}\\
\,u(-a)=1,\qquad u(0)=\theta,\qquad  u(a)=0,\vspace{3pt}\\
u>0 \quad\text{ in } (-a,a),
\end{cases}
$$
which is such that
$$
\Vert u \Vert _{C^2(-a,a)}\leq K, \quad -c_{min}\leq c \leq
c_{max}.
$$
\end{prop}

\noindent {\bf Proof.} For a given nonnegative function $v$
defined on $(-a,a)$ and satisfying the Dirichlet boundary
conditions as requested in $P_\tau(a)$ --- that is $v(-a)=1$ and
$v(a)=0$--- consider the family $0\leq\tau\leq1$ of linear
problems
\begin{equation}\label{droite-gelee}
P_\tau^c(a) \;\begin{cases}\, -u''-cu'=\tau v(v-\theta)(1-\phi * v) \quad &\text{ in } (-a,a)\vspace{3pt} \\
u(-a)=1,\quad u(a)=0.
\end{cases}
\end{equation}
Denote by $\mathcal K _\tau$ the mapping of the Banach space
$X:=\R \times C^{1,\alpha}(Q)$ --- equipped with the norm $\Vert
(c,v)\Vert_X:=\max\left(|c|,\Vert v\Vert
_{C^{1,\alpha}}\right)$--- onto itself defined by
$$
\mathcal K _\tau:(c,v)\mapsto \left(\theta-v(0)+c,u_\tau
^c:=\text{ the solution of } P_\tau ^c(a)\right).
$$
Constructing a solution $(c,u)$ of $P_1(a)$ is equivalent to
showing that the kernel of $\text{Id} -\mathcal K _1$ is
nontrivial. The operator $\mathcal K _\tau$ is compact and depends
continuously on the parameter $0\leq \tau \leq 1$. Thus the
Leray-Schauder topological argument can be applied. Define the
open set
$$
S:=\left\{(c,v):\, -\tilde c_{min}(a)-1< c< c_{max}+1,\;v>0,\;
\Vert v \Vert _{C^{1,\alpha}}< M+1\right\}\subset X.
$$
It follows from the {\it a priori} estimates Lemma
\ref{lem:a-priori-bounds}, Lemma \ref{lem:a-priori-speed} and
Lemma \ref{lem:a-priori-speed3}, that there exists $a_0>0$ such
that, for any $a\geq a_0$, any $0\leq \tau \leq 1$, the operator
$\text{Id} -\mathcal K _\tau$ cannot vanish on the boundary
$\partial S$. By the homotopy invariance of the degree we thus
have $\text{deg}(\text{Id}-\mathcal
K_1,S,0)=\text{deg}(\text{Id}-\mathcal K_0,S,0)$.

To conclude, observe that we can compute
\begin{equation}\label{tauzero}
u_0^c(x)=\displaystyle \frac{e^{-cx}-e^{-ca}}{e^{ca}-e^{-ca}}
\quad\text{ if } c\neq 0,\quad u_0^c(x)=-\frac 1{2a}x+\frac 12
\quad\text{ if } c=0,
\end{equation}
and that $u_0^c(0)$ is decreasing with respect to $c$ (in
particular there is a unique $c_0$ such that
$u_0^{c_0}(0)=\theta$). Hence by using two additional homotopies
(see \cite{Ber-Nad-Per-Ryz} or \cite{A-Cov-Rao} for details) we
can compute
 $\text{deg}(\text{Id}-\mathcal
K_0, S,0)=-1$
 so that $\text{deg}(\text{Id}-\mathcal
K_1,S,0)=-1$ and there is a
 $(c,u)\in S$ solution of $P_1(a)$. Finally, Lemma \ref{lem:a-priori-speed2} provides a lower bound
 $c\geq -c_{min}$, uniform in $a\geq a_0$. This concludes the proof of
the proposition. \fdem

 \medskip

 \noindent{\bf A solution on $\R$.} Equipped with the solution $(c,u)$ of $P_1(a)$ of Proposition
\ref{prop:sol-boite}, we now let $a\to +\infty$. This enables to
construct
--- passing to a subsequence $a_n \to +\infty$--- a speed
$-c_{min}\leq \c \leq c_{max} $ and a function $U:\R\to(0,M)$ in
$C^2 _b (\R)$ such that
\begin{equation}\label{eq-onde-construite}
-U''-\c U^{\prime}=U(U-\theta)(1-\phi *  U)\quad \text{ on }\R,
\end{equation}
\begin{equation}\label{masse-onde-construite}
U(0)=\theta.
\end{equation}

In contrast with the nonlocal Fisher-KPP equation considered in
\cite{Ber-Nad-Per-Ryz} we need additional arguments to show that
the constructed $U$ is non trivial, i.e. that $U$ \lq\lq
visits\rq\rq\, both $(0,\theta)$ and $(\theta,1)$. This is the
purpose of the next section.

\section{Non triviality of $U$ the solution on $\R$ }\label{s:nontrivial}

In this section, we provide additional {\it a priori} estimates on
the solution $u$ in the box $(-a,a)$, which in turn will imply the
non triviality of the solution $U$ on $\R$.

\medskip
First, using the homotopy of the previous section, we show that
the solution in the box cannot attain $\theta$ elsewhere that at
$x=0$.

\begin{prop}[$\theta$ is attained only at $x=0$]\label{u0iffx0}
For all $a\geq a_0$, the solution $(c,u)$ of Proposition
\ref{prop:sol-boite} satisfies
$$
u(x)=\theta \; \text{ if and only if }\;  x=0.
$$
\end{prop}

\noindent {\bf Proof.} From Proposition \ref{prop:sol-boite}, we
know that there is a solution $(c_\tau,u_\tau)$ of
\begin{equation}\label{twtau}
\left\{\begin{array}{l}
-u_\tau''-c_\tau u_\tau'=\tau u_\tau(u_\tau-\theta)(1-\phi *  {u_\tau})\quad \text{ in }(-a,a)\vspace{3pt}\\
 u_\tau(-a)=1,\qquad u_\tau(0)=\theta,\qquad  u_\tau(a)=0,
  \end{array}
\right.
\end{equation}
and that $(c_\tau,u_\tau)$ depends continuously upon $0\leq \tau
\leq 1$. For $\tau =0$, in view of \eqref{tauzero}, the solution
$u_0$ satisfies $u_0(x)=\theta$ if and only if $x=0$. We can
therefore define
$$
\tau^*:=\sup \left\{ 0\leq \tau \leq 1, \forall \sigma \in
[0,\tau], u_\sigma(x)=\theta \; \text{ iff } \; x=0 \right\}.
$$

 Assume by contradiction that there is a
$x^*\neq 0$ such that $u_{\tau ^*}(x^*)=\theta$. Without loss of
generality, we can assume $x^*<0$ and $u_{\tau ^*}>\theta$ on
$(x^*,0)$. By the definition of $\tau ^*$ as a supremum, one must
have $u_{\tau ^*}\geq \theta$ on $(-a,0)$, which in turn enforces
$u'_{\tau^*}(x^*)=0$. Hence $v:=u_{\tau ^*}-\theta$ is positive on
$(x^*,0)$, zero at $x^*$ and satisfies the linear elliptic
equation
$$
-v''-c_{\tau ^*}v'=\left[\tau ^*u_{\tau^*}(1-\phi*
u_{\tau^*})\right]v\quad \text{ on } (x^*,0).
$$
It then follows --- see e.g. \cite[Lemma 3.4]{Gil-Tru}--- that
$v'(x^*)>0$, which is a contradiction. Hence $u_{\tau^*}$ attains
$\theta$ only at $x=0$. To conclude let us prove that $\tau ^*=1$.

Assume by contradiction that $0\leq \tau ^*<1$. By the definition
of $\tau ^*$, there exists a sequence $(\tau_n,x_n)$ such that
$\tau_n\downarrow \tau^\ast$, $x_n\neq 0$, and
$u_{\tau_n}(x_n)=0$. Up to an extraction, the sequence $x_n$
converges to a limit $x^\ast$, which implies, thanks to the
continuity of $(\tau,x)\mapsto u_\tau(x)$ with respect to $\tau$
and $x$, that $u_{\tau^\ast}(x^\ast)=0$. As seen above one must
have $x^*=0$, and then $x_n\to 0$. Then, for some $-1\leq C_n \leq
1$, we have
$0=u_{\tau_n}(x_n)=u_{\tau_n}(0)+u_{\tau_n}'(0)x_n+C_n\|u_{\tau_n}\|_{C^{1,\alpha}}|x_n|^{1+\alpha}$,
that is
$$
|u_{\tau_n}'(0)|\leq
|C_n|\,\|u_{\tau_n}\|_{C^{1,\alpha}}|x_n|^{\alpha}\leq C
|x_n|^\alpha \to_{n\to\infty}0.
$$
The continuity of $(u_\tau')$ with respect to $\tau$ then implies
that $u_{\tau ^*}'(0)=0$. Since $u_{\tau ^*}>\theta$ on $(-a,0)$
we derive a contradiction as above. As a result $\tau ^*=1$ and
the proposition is proved.\fdem

\medskip
We now construct a subsolution of a linear equation, having the
form of a bump, that will be very useful in the following.

\begin{lem}[A  bump as a sub-solution]\label{functionphi} For any $\kappa >0$, there exists $A>0$ such that for all
$c>-2\sqrt\kappa$, there exist $0<\tilde x<A$, $X>\tilde x$ and
$\psi:[0,X]\to [0,1]$, satisfying $\psi(0)=0$, $\psi(\tilde x)=1$
and
 \begin{equation}\label{eqphi}
-\psi ''-c \psi'\leq \kappa \psi \quad \text{ in }(0,X).
 \end{equation}
\end{lem}

 \noindent {\bf Proof.} If $-2\sqrt\kappa< c <2\sqrt\kappa$, we define $ \psi(x):=e^{-\frac c 2  x}
 \sin \left(\frac {\sqrt{4\kappa -c^2}}2x\right)$ which solves
 $- \psi ''-c \psi '=\kappa  \psi$. Also on
 $[0,2\pi/\sqrt{4\kappa-c^2}]$ we have
 $ \psi(0)=\psi(2\pi/\sqrt{4\kappa-c^2})=0$, $
 \psi \geq 0$ and maximal at point
 $$
 \tilde x=\tilde x(c)=\begin{cases}\frac 2{\sqrt{4\kappa
 -c^2}}\tan ^{-1} \left(\frac{\sqrt{4\kappa -c^2}}c\right)& \text{ if } 0<c<2\sqrt \kappa \\
 \frac \pi {\sqrt{4\kappa}} & \text{ if } c=0\\
\frac 2{\sqrt{4\kappa
 -c^2}}\left(\tan ^{-1} \left(\frac{\sqrt{4\kappa
 -c^2}}c\right)+\pi\right)& \text{ if } -2\sqrt \kappa <c<0.
 \end{cases}
 $$
$\tilde x(c)$ is then uniformly bounded for $c\in(-2\sqrt
\kappa,2\sqrt \kappa)$, so that the renormalized function
$\psi/\psi(\tilde x)$ is as requested.

If $c\geq 2 \sqrt \kappa$, we define
$\psi(x):=e^{\frac{-\sqrt\kappa}2x}\sin\left(\frac{\sqrt\kappa}2x\right)$
and $\tilde x=\pi/\sqrt{4\kappa}$, so that $\psi$ increases on
$[0,\tilde x]$ and starts to decreases after $\tilde x$. On
$[0,\tilde x]$, we have $\tan
\left(\frac{\sqrt\kappa}2x\right)\leq 1$ so that
\begin{eqnarray*}
 -\psi''(x)-c\psi'(x)-\kappa\psi(x)&=&\sqrt{\kappa}e^{\frac{-\sqrt\kappa}2x}\cos
 \left(\frac{\sqrt\kappa}2x\right)
 \left(\left(\frac{\sqrt\kappa}2-\frac{c}2\right)+\left(\frac{c}2-\sqrt\kappa\right)\tan\left(\frac{\sqrt\kappa}2x\right)\right)\\
 &\leq & \frac{\sqrt \kappa} 2e^{\frac{-\sqrt\kappa}2x}\cos
 \left(\frac{\sqrt\kappa}2x\right)(\sqrt \kappa -c +c-2\sqrt
 \kappa)\\
 &\leq& 0.
 \end{eqnarray*}
Observe that $-\psi''(\tilde x)-c\psi'(\tilde x)-\kappa\psi(\tilde
x)\leq -\frac {\kappa e^{-\pi/4}\cos(\pi /4)} 2<0$, so that there
is $X>\tilde x$ such that $-\psi ''-c\psi '-\kappa \psi \leq 0$ on
$[0,X]$. Hence $\psi / \psi(\tilde x)$ is as requested. \fdem

\medskip

We will also use the elementary following lemma.

\begin{lem}[An auxiliary solution]\label{functionchi} Let $\rho>0$ and $b>0$ be given.
Then, for all $c\in \R$,  there is a decreasing function
$\chi=\chi _c:\R\to \R$ such that $\chi(0)=1$, $\chi(b)=0$ and
 \begin{equation}\label{eqchi}
-\chi ''-c \chi' =-\rho \chi \quad \text{ in } \R.
 \end{equation}
\end{lem}

\noindent {\bf Proof.} One solves the linear ODE and sees that the
function
$$
\chi(x):=\left(1-\frac
1{1-e^{-\sqrt{c^2+4\rho}\,b}}\right)e^{\frac{-c+\sqrt{c^2+4\rho}}2x}+\frac
1{1-e^{-\sqrt{c^2+4\rho}\,b}}e^{\frac{-c-\sqrt{c^2+4\rho}}2x}
$$
is as requested. \fdem

\medskip We now show that $u$ can be uniformly (with respect to $a$) bounded away from $\theta$ far on the
right or the left, depending on the sign of the speed $c$.

\begin{prop}[Moving away from $\theta$] \label{propboundaway}
There exist $\ep >0$ and $a_0>0$ such that, for all $a\geq a_0$,
any solution $(c,u)$ of
\begin{equation}\label{equ}
\begin{cases}
-u''-cu^{\prime}=u(u-\theta)(1-\phi *  u)\quad \text{ in }(-a,a)\vspace{3pt}\\
u(-a)=1,\qquad u(0)=\theta,\qquad  u(a)=0,\vspace{3pt}\\
u(x)=\theta \; \text{ if and only if }\;  x=0,
\end{cases}
\end{equation}
satisfies, if we define $\kappa:=\theta(1-\theta)/8 >0$,
\begin{itemize}
\item [(i)] $c>-2\sqrt\kappa \Longrightarrow u\leq \theta /2$ on
$[1/\ep,a]$ \item [(ii)] $c<2\sqrt\kappa \Longrightarrow u\geq
\theta+\ep$ on $[-a,-1/\ep]$.
\end{itemize}
\end{prop}

\noindent {\bf Proof.} For $\kappa =\theta (1-\theta)/8$, let
$\psi$ be the bump of Lemma \ref{functionphi}.

\medskip

Assume $c>-2\sqrt\kappa$ and let us prove $(i)$.  Since $\phi\in
L^1$, there exists $R>0$ such that $\int_R^\infty\phi\leq
\frac{1-\theta}{2M}$, where $M$ is the $L^\infty$ bound we have on
$u$. We turn upside down the bump and make it slide from the right
towards the left until touching the solution $u$. Precisely one
can define
$$
\alpha_0:=\min \left\{\alpha\geq 0:\, \forall x-\alpha \in[0,X],\,
u(x)<  \theta -\frac\theta 2\psi(x-\alpha)\right\}\in [0,a).
$$
We aim at proving that $\alpha _0 \leq R$ uniformly with respect
to large $a$. Assume by contradiction that $\alpha _0 >R$. The
function $v:=u-\theta+\frac \theta 2 \psi (\cdot -\alpha _0)$ has
a zero maximum at some point $x_0$. Notice that since $\psi(\tilde
x)=\max \psi$, the definition of $\alpha_0$ implies that
$x_0-\alpha_0\in (0,\tilde x]$, so that $\psi$ is a subsolution of
\eqref{eqphi} around $x_0-\alpha_0$. Thus $0\geq
v''(x_0)+cv'(x_0)$ implies
\begin{eqnarray}
 0&\geq& (u''+cu')(x_0)+\frac \theta 2(\psi''+c\psi')(x_0-\alpha _0)\nonumber\\
 &\geq& -u(x_0)(u(x_0)-\theta)(1-\phi\ast  u(x_0))-\frac {\kappa\theta} 2\psi(x_0-\alpha_0)\nonumber\\
&\geq&(\theta-u(x_0))\left[u(x_0)(1-\phi\ast
u(x_0))-\kappa\right].\label{truc}
\end{eqnarray}
Now, since $\alpha_0\geq R$, we have $x_0\geq R$ and we can
estimate the nonlocal term by
\begin{eqnarray}
 \phi\ast u(x_0)&\leq&\int_{-\infty}^0\phi(x_0-y) u(y)\,dy+\int_0^{\infty}\phi( x_0-y) u(y)\,dy\nonumber\\
 &\leq&M\int_R^\infty\phi+\theta \int_{\R}\phi\leq\frac{1+\theta}{2}.\label{estimconv1}
\end{eqnarray}
Since $\frac \theta 2 \leq u(x_0)<\theta$, it follows from
\eqref{truc} that $0\geq \frac \theta 2 \frac{1-\theta }2
-\kappa$, which contradicts the definition of $\kappa$. As a
result $\alpha _0 \leq R$, which means that the minimum $\theta
/2$ of the reversed bump can slide to the left at least until
$R+\tilde x$. In other words we have $u\leq \theta /2$ on
$[R+\tilde x,a]$ which concludes the proof of $(i)$.

\medskip
Assume $c<2\sqrt\kappa$ and let us prove $(ii)$. Since $\phi\in
L^1$, we can choose $R>0$ such that $\int_{[-R,R]^c}\phi\leq
\frac{1-\theta}{4M}$, where $M$ is the $L^\infty$ bound we have on
$u$. Before using the bump we need a preliminary result via the
function $\chi$ of Lemma \ref{functionchi}.

For $\rho:=2M$ and $b:=2R+1$ define $\chi$ as in Lemma
\ref{functionchi}. Provided that $a>2b$, Proposition \ref{u0iffx0}
shows that, for $\lambda >0$ small enough, $\theta+\lambda
\chi(\cdot+a)<u$ on $[-a,-a+b]$. We can therefore define
$$
\lambda_0:=\max\left\{\lambda>0:\,\forall x\in[-a,-a+b],\,
\theta+\lambda\, \chi(x+a) < u(x)\right\}\in (0,1-\theta].
$$
The function $v:=u-\theta-\lambda_0\chi(\cdot +a)$ thus has a zero
minimum at a point $x_0$. Assume by contradiction that $\lambda _0
<1-\theta$, which in turn implies $x_0\neq -a$. Also Proposition
\ref{u0iffx0} implies $u(b)>\theta$ so that $x_0 \neq -a+b$. Thus
$0\leq v''(x_0)+cv'(x_0)$ so that
\begin{eqnarray*}
0&\leq&(u''+cu')(x_0)-\lambda_0\,(\chi''+c\chi')(x_0+a)\\
&=&-u(x_0)(\theta-u(x_0))(1-\phi\ast  u(x_0))-\lambda_0\,\rho\,\chi(x_0+a)\\
&=&(u(x_0)-\theta)\left[u(x_0)(1-\phi\ast
u(x_0))-\rho\right]\\
&\leq &(u(x_0)-\theta)(M-\rho)<0,
\end{eqnarray*}
which is absurd. Hence $\lambda _0=1-\theta$ and thus
\begin{equation}\label{touchepasaudebut}
u(x)\geq \theta+(1-\theta)\chi(x+a),\, \forall x\in
[-a,-a+b]=[-a,-a+2R+1].
\end{equation}

Let us now define, for $\ep >0$ to be selected,
$$
\alpha_0:=\max \left\{\alpha\leq 0:\, \forall \alpha -x
\in[0,X],\, u(x)>  \theta +\ep\psi(\alpha-x)\right\}\in (-a,0].
$$
The estimate \eqref{touchepasaudebut} shows that it is enough to
choose $\ep<(1-\theta)\min_{[-a,-a+2R]}\chi(\cdot+a)$ to get the
lower bound $\alpha_0 \geq -a+2R$. We aim at proving that $\alpha
_0 \geq -2R$ uniformly with respect to large $a$. Assume by
contradiction that $\alpha _0 <-2R$. The function $v:=u-\theta-\ep
\psi (\alpha _0-\cdot )$ has a zero minimum at some point $x_0$.
Notice that since $\psi(\tilde x)=\max \psi$, the definition of
$\alpha_0$ implies that $\alpha_0-x_0 \in (0,\tilde x]$, so that
$\psi$ is a subsolution of \eqref{eqphi} around $\alpha_0-x_0$.
Thus, we have
\begin{eqnarray}
 0&\leq& (u''+cu')(x_0)-\ep(\psi''+c\psi')(\alpha _0 -x_0)\nonumber\\
 &\leq& -u(x_0)(u(x_0)-\theta)(1-\phi\ast  u(x_0))+\ep\kappa \psi(\alpha_0-x_0)\nonumber\\
&\leq&(u(x_0)-\theta)\left[\kappa-u(x_0)(1-\phi\ast
u(x_0))\right].\label{bidule}
\end{eqnarray}
Now observe that $-a+2R \leq \alpha _0 <-2R$ implies $-a+2R \leq
x_0 \leq -2R$, so that $[x_0-R,x_0+R]\subset [-a+R,-R]$. Therefore
the Harnack inequality applied to $u-\theta$ provides a constant
$C>0$, independent of $a$, such that
$$
0<u(x)-\theta \leq C(u( x_0)-\theta),\, \forall x\in
[x_0-R,x_0+R].
$$
This allows to estimate the nonlocal term by
\begin{eqnarray}
 \phi\ast  u(x_0)&\leq&\int_{[-R,R]}\phi(y) u(x_0-y)\,dy+\int_{[-R,R]^c}\phi(y) u(x_0-y)\,dy\nonumber\\
 &\leq& \theta+C(u(x_0)-\theta)+\frac{1-\theta}{4}\leq \frac{1+\theta}{2},\label{estimconv2}
\end{eqnarray}
provided that $u(x_0)\leq \theta+\frac{1-\theta}{4C}$, which is
satisfied if we choose $\ep>0$ small enough (we recall that
$u(x_0)\leq \theta+\ep\|\psi\|_\infty=\theta +\ep$). It follows
that $\kappa-u(x_0)(1-\phi\ast  u(x_0))\leq
-\frac{\theta(1-\theta)}4<0$, which contradicts \eqref{bidule}. As
a result $\alpha_0\geq -2R$, which concludes the proof. \fdem

\medskip
\noindent{\bf Non triviality of $U$.} Let us recall that $(\c,U)$
is constructed as the limit of $(c_a,u_a)$ as $a\to \infty$. By
extraction if necessary we can assume that
 the $(c_a,u_a)$' satisfy either $(i)$ or $(ii)$ of Proposition
 \ref{propboundaway}, and so does $(\c,U)$. As a result, the
 constructed wave $(\c,U)$ is non trivial.

\section{Behaviors of $U$ in $\pm \infty$}\label{s:properties}

We now prove the behavior \eqref{gauche} as $x\to -\infty$, the
limit \eqref{droite} as $x\to \infty$, and that the constructed
front is decreasing for $x>0$ large enough. This will complete the
proof of Theorem \ref{th:tw}.

\begin{prop}[Behaviors of $U$ at infinity]\label{prop:uinfini}
 Let $(\c,U)$ be the solution of \eqref{eq-dans-espace} constructed in the end of Section \ref{s:tw}.
 Then, for some $\ep >0$,
$$
U(x)\geq \theta +\ep\; \text{ for all } x\in(-\infty,-1/\ep),
$$
 and
$$
\lim_{x\to+\infty}U(x)=0.
$$
Moreover, there exists $\bar x>0$ such that $U$ is decreasing on
$[\bar x,\infty)$.
\end{prop}

\noindent {\bf Proof.} \emph{Step 1: We show that $U > \theta$ on
$(-\infty,0)$, and $U<\theta$ on $(0,\infty)$.}

 Thanks to Lemma \ref{u0iffx0}, for
any $a>0$, the solution $(c,u)$ in the box of Proposition
\ref{prop:sol-boite} satisfies $u\geq\theta$ on $[-a,0]$, and
$u\leq\theta$ on $[0,a]$. Since $(\c,U)$ is, on any compact
interval, the uniform limit of such solutions, it satisfies $U
\geq \theta$ on $(-\infty,0]$ and $U\leq\theta$ on $[0,\infty)$.
Thus, any $x\neq 0$ such that $U(x)=\theta$ is a local extremum
and $U'(x)=0$; this is impossible, since  $U$ is a solution of
\eqref{eq-dans-espace} and $U\not\equiv\theta$ thanks to Section
\ref{s:nontrivial}.

\medskip

\noindent\emph{Step 2: We show that there exists $\ep>0$ such that
$|U-\theta|\geq\ep$ on $[-1/\ep,1/\ep]^c$.}

 Consider first the case where $\c\geq 0$. Since $U$ is a
limit of solutions of Proposition \ref{prop:sol-boite}, the
Proposition \ref{propboundaway} shows that there exists $\ep>0$
such that $U\leq\theta /2$ on $(1/\ep,\infty)$. To investigate the
left side, as in the proof of Proposition \ref{propboundaway}
$(ii)$, we choose $R>0$ such that $\int_{[-R,R]^c}
\phi\leq\frac{1-\theta}{4M}$,
 where $M$ is the $L^\infty$ bound we have on $U$. By the Harnack inequality applied to $U-\theta$, there exists $C>0$ such
 that, for all $x_0\leq -2R$,
\begin{equation}\label{harnack}
0<U(x)-\theta \leq C(U( x_0)-\theta),\, \forall x\in
[x_0-R,x_0+R].
\end{equation}
Let us now define, for $\eta \geq 0$,
$$
\psi_\eta(x):=\theta+\gamma(1-\eta (-2R-x)),\quad \gamma:=
\min\left(\gamma_1:=\frac
12(U(-2R)-\theta),\gamma_2:=\frac{1-\theta}{4C}\right)>0.
$$
Since $\psi_\eta\leq U$ on $(-\infty,-2R]$ for $\eta>0$ large
enough,  we can define
$$\eta_0:=\min\left\{\eta\geq 0:\, \forall x\leq -2R,\, \psi_\eta(x)\leq U(x)\right\}.$$
Let us assume by contradiction that $\eta_0>0$. The function
$U-\psi_{\eta_0}$ then attains a zero minimum at a point $x_0<
-2R$  (notice that $\gamma \leq \gamma _1$ prevents $x_0=-2R$).
Hence
\begin{eqnarray*}
 0&\geq& -(U-\psi_{\eta_0})''(x_0)-\c(U-\psi_{\eta_0})'(x_0)\\
 &\geq&\c\psi_{\eta_0}'(x_0)+U(x_0)(U(x_0)-\theta)(1-\phi\ast U(x_0))\\
 &\geq&\c\gamma \eta_0+U(x_0)(U(x_0)-\theta)\frac{1-\theta}2>0,
\end{eqnarray*}
where we have used the estimate \eqref{estimconv2} for $U$ (notice
that this is possible since we have the two ingredients
\eqref{harnack} and $U(x_0)=\psi_{\eta _0}(x_0)\leq \theta
+\frac{1-\theta}{4C}$). This is a contradiction which proves that
$\eta_0=0$, and then $U\geq \theta +\gamma$ on $(-\infty,-2R)$.
This concludes the case $\c \geq 0$.

\medskip

Consider next the case where $\c\leq 0$. Since $U$ is a limit of
solutions of Proposition \ref{prop:sol-boite}, the Proposition
 \ref{propboundaway} shows that there exists $\ep>0$ such that $U\geq \theta+\ep$ on $(-\infty,1/\ep)$. To
 investigate the right side, as in the proof of
 Proposition \ref{propboundaway} $(i)$, we choose $R>0$ such that $\int_R^\infty \phi\leq\frac{1-\theta}{2M}$, were $M$ is the $L^\infty$
 bound we have on $U$. We define
 $$
 \psi_\eta(x):=\theta+\gamma(-1+\eta (x-2R)),\quad \gamma:=\frac
 12 (\theta -U(2R))>0,
 $$ which satisfies $\psi_\eta\geq U$
  on $[2R,\infty)$ for $\eta>0$ large enough. We can then define
$$\eta_0:=\min\left\{\eta\geq 0:\forall x\geq 2R,\, \psi_\eta(x)\geq U(x)\right\}.$$
Let us assume by contradiction that $\eta_0>0$. The function
$U-\psi_{\eta_0}$ then attains a zero maximum at a point $x_0>2R$,
 and therefore
\begin{eqnarray*}
 0&\leq& -(U-\psi_{\eta_0})''(x_0)-\c(U-\psi_{\eta_0})'(x_0)\\
 &\leq&\c\psi_{\eta_0}'(x_0)+U(x_0)(U(x_0)-\theta)(1-\phi\ast U(x_0))\\
 &<&\c\gamma \eta_0+U(x_0)(U(x_0)-\theta)\frac{1-\theta}2<0,
\end{eqnarray*}
where we have used the estimate \eqref{estimconv1} for $U$. This
is a contradiction which proves that $\eta_0=0$, and then $U\leq
\theta-\gamma$ on $(2R,\infty)$. This concludes the case $\c\leq
0$.

\medskip

\noindent\emph{Step 3: We show that $U$ decreases to $0$ on some
interval $(\bar x,\infty)$.}

Choose $R>0$ large enough so that
$\int_R^\infty\phi\leq\frac{1-\theta}{2M}$. Assume by
contradiction that $U$ admits a local maximum at some point $x_m
\geq R$. Since $U(x_m)<\theta$, by evaluating the equation
\eqref{eq-dans-espace} we see that $1\leq \phi * U (x_m)$. But on
the other hand
$$
\phi * U(x_m) \leq \int _{-\infty}^0 \phi (x_m-y)u(y)\,dy+\int _0
^\infty \phi(x_m-y)u(y)\,dy\leq M \int _R ^\infty \phi + \theta
\leq \frac{1+\theta}2<1,
$$
which is a contradiction. Hence $U$ cannot attain a maximum on
$(R,\infty)$, which in turn implies that there is $\bar x >0$ such
that $U$ is monotonic (increasing or decreasing) on $[\bar
x,\infty)$. Hence, as $x\to \infty$, $U(x)\to l$ and, by Step 2,
$0\leq l\leq \theta -\ep$.

We define next $v_n(x):=U(x+n)$, which solves
$$
-{v_n}''-\c{v_n}'=v_n(v_n-\theta)(1-\phi*v_n)\quad \textrm{ on }
\R.
$$
Since the $L^\infty$ norm of the right hand side member is
uniformly bounded with respect to $n$, the interior elliptic
estimates imply that, for all $R>0$, all $1<p<\infty$, the
sequence $(v_n)$ is bounded in $W^{2,p}([-R,R])$. From Sobolev
embedding theorem, one can extract $v_{\varphi (n)} \to v$
strongly in $C^{1,\beta}_{loc}(\R)$ and weakly in $W^{2,p}_{loc}
(\R)$. Since $v_n(x)=U(x+n)\to l$ we have $v\equiv l$ and
$v'\equiv 0$. Combining this with the fact that $v$ solves
$$
-v''-\c v'= v(v-\theta)(1-\phi *v)\quad \textrm{ on } \R,
$$
we have $l(l-\theta)(1-l)=0$, which implies $l=0$ and the decrease
of $U$ on $[\bar x,\infty)$.\fdem

\section{Focusing kernels}\label{s:focusing}

In this section we consider $(\c _\sigma, U_\sigma)$ the
constructed waves for the focusing kernels
$$
\phi _\sigma (x)=\frac 1 \sigma \phi \left(\frac x \sigma\right),
\quad \sigma>0.
$$
We prove Proposition \ref{prop:focusing}. Item $(i)$ consists in a
perturbation analysis, and item $(ii)$ will follow from the $L^2$
analysis performed in \cite{Alf-Cov}.

\medskip

\noindent{\bf Proof of $(i)$.} Assume $m_1:=\int _\R
|z|\phi(z)\,dz<\infty$. We have
\begin{equation}\label{eq-dans-espace_sigma}
-{U_\sigma}''-\c _\sigma {U_\sigma}'=
U_\sigma(U_\sigma-\theta)(1-\phi _\sigma
* U_\sigma) \quad \text{ on } \R,
\end{equation}
and $0\leq U_\sigma \leq M_\sigma$, $c_{min,\sigma} \leq \c
_\sigma \leq c_{max,\sigma}$, with $M_\sigma$, $c_{min,\sigma}$, $
c_{max,\sigma}$ depending {\it a priori} on $\sigma
>0$. The following lemma improves the bounds for the
travelling waves: as $\sigma \to 0$ solutions
$(\c_\sigma,U_\sigma)$ are uniformly bounded.

\begin{lem}[Uniform bounds for $(\c_\sigma,U_\sigma)$]\label{lem:a-priori-bounds-focusing}
 Let $\sigma _0>0$ be arbitrary. Then there is $M>0$, $ c_{min} \in \R$, $ c_{max}\in \R$ such
  that, for all $\sigma\in(0,\sigma _0)$,
$$
0\leq U_\sigma \leq  M,\quad\text{ and } \quad  c_{min}\leq
\c_\sigma \leq  c_{max}.
$$
\end{lem}

\noindent {\bf Proof.} It is sufficient to work on the solutions
$(c_\sigma,u_\sigma)$ in the box. Define
$M_\sigma:=\max_{x\in[-a,a]} u_\sigma(x)$. A first lecture of
Lemma \ref{lem:a-priori-bounds} yields the rough bound
\eqref{borne} with the kernel $\phi _\sigma$ in place of $\phi$.
Since $\Vert \phi _\sigma '\Vert_{L^\infty(-1,1)}\leq \frac 1
{\sigma ^2}\Vert \phi '\Vert _{L^\infty(\R)}$ and $\phi _\sigma
(0)=\frac 1 \sigma \phi(0)$, we infer from \eqref{borne} that
there is a constant $b
>0$, such that $M _\sigma\leq
b ^2/\sigma ^2$. Equipped with this rough bound, we go back to the
proof of Lemma \ref{lem:a-priori-bounds} but rather than
\eqref{choix} we select the improvement $$ x_0:=\frac 1{2b}
\sigma.
$$
Hence, going further into the proof, we discover
$$
1\geq \frac 3 4M_\sigma \int _0 ^{\frac{\sigma}{2b}} \phi
_\sigma=\frac 34 M_\sigma \int _0 ^{\frac 1 {2b}} \phi,
$$
so that $M_\sigma \leq  M:=\frac 4 3 \left( \int _0 ^{\frac 1{2b}}
\phi\right)^{-1}$, that is a uniform bound $M$ for $M_\sigma$ as
$\sigma \to 0$.

In view of Lemma \ref{lem:a-priori-speed} and of the proof of
Lemma \ref{lem:a-priori-speed2}, the uniform bound $M$ yields
uniform bounds $c_{max}$ and $c_{min}$ for the speed $c_\sigma$.
The lemma is proved. \fdem

\medskip

Hence, the coefficients and the right hand side member of the
elliptic equation \eqref{eq-dans-espace_sigma} are uniformly
bounded w.r.t. $\sigma \in(0,\sigma _0)$. Therefore Schauder's
elliptic estimates
--- see \cite[(1.11)]{Lad-Ura} for instance--- imply that, $\Vert U_\sigma\Vert
_{C^{2,\beta}}\leq C_0$ with $C_0>0$ not depending on $\sigma$. It
follows that
$$
|U_\sigma-\phi _\sigma *U_\sigma|(x)\leq \int _{\R} \phi _\sigma
(x-y)|U_\sigma(x)-U_\sigma (y)|\,dy \leq \Vert {U_\sigma}'\Vert
_{\infty}\int _{\R}\phi_\sigma(x-y)|x-y|\,dy \leq C_0 m_1 \sigma.
$$
Hence writing $1-\phi _\sigma * U_\sigma=1-U_\sigma+U_\sigma
-\phi_\sigma * U_\sigma$ in \eqref{eq-dans-espace_sigma}, we get,
for some $C>0$,
$$
f^+_\sigma(U_\sigma) \geq -{U_\sigma}''-\c _\sigma {U_\sigma}'\geq
f^-_\sigma(U_\sigma) \quad \text{ on } \R,
$$
where
$$
f^\pm _\sigma (s):=s(s-\theta)(1-s)\pm C\sigma.
$$
Hence, by the comparison principle, $\psi ^+_\sigma(x,t)\geq
U_\sigma (x-\c _\sigma t)\geq \psi ^- _\sigma (x,t)$, with $\psi
^\pm _\sigma$ the solutions of the Cauchy parabolic problems
\begin{equation*}
\begin{cases}
\partial _t \psi=\partial _{xx} \psi +f^\pm _\sigma(\psi)\quad \text{ in }(0,\infty)\times \R,\\
\psi(x,0)=U_\sigma (x)\quad\text{ in } \R.
\end{cases}
\end{equation*}
Observe that, for $\sigma >0$ small enough, the functions $f^\pm
_\sigma$ are still of the bistable type with three zeros
$\alpha^\pm _\sigma=\mathcal O(\sigma)$, $\beta ^\pm _\sigma
=\theta+\mathcal O(\sigma)$, $\gamma ^\pm _\sigma=1+\mathcal O
(\sigma)$. It is therefore well-known \cite[Theorem 3.1]{Fif-Mac}
that, for a given small $\sigma>0$, the solutions $\psi ^\pm
_\sigma$ approach $U^\pm _\sigma(x-c^\pm _\sigma t-x_0 ^\pm)$, for
two given $x_0^\pm \in \R$, uniformly in $x$ as $t\to \infty$.
Here $(c^\pm _\sigma,U^\pm _\sigma)$ denotes the bistable wave
\begin{equation*}
\begin{cases}
{U^\pm _\sigma}''+c^\pm _\sigma{U^\pm _\sigma}'+f^\pm _\sigma (U^\pm _\sigma)=0,\\
\lim_{x\to-\infty}U^\pm _\sigma(x)=\gamma ^\pm _\sigma, \quad
U^\pm _\sigma(0)=\beta ^\pm _\sigma, \quad \lim_{x\to+\infty}U^\pm
_\sigma (x)=\alpha ^\pm _\sigma.
\end{cases}
\end{equation*}
This enforces
$$
c^+_\sigma \geq \c _\sigma \geq c^- _\sigma.
$$
Since, as $\sigma \to 0$, $c_\sigma ^\pm$ converge to $\c _0$ the
speed of the wave \eqref{local}, this concludes the proof of
$(i)$.

\medskip

\noindent{\bf Proof of $(ii)$.}  Assume $\theta \neq \frac 12$,
which in turn implies $\c _0\neq 0$, and $m_2:=\int _\R z^2
\phi(z)\,dz<\infty$. Observe that $\int _\R z^2 \phi _\sigma
(z)\,dz=\sigma ^2 m_2$ so that, in virtue of \cite[Lemma
5]{Alf-Cov}, to get $\lim_{x\to-\infty}U_\sigma (x)=1$ it is
enough to have
\begin{equation}\label{enough}
\sigma \sqrt m_2 {M_\sigma}^2 < |\c_\sigma|,
\end{equation}
which is clear, for small enough $\sigma>0$, since $M_\sigma \leq
M$, and $|\c _\sigma|\to |\c _0|\neq 0$. \fdem

\section{The ignition case}\label{s:ignition}

Here we explain briefly how to use  similar arguments to handle
the case of the ignition case.

\medskip

The typical local ignition case is given by $-U''-\c U'=\mathbf 1
_{\{U\geq \theta\}} (U-\theta)(1-U)$, and the corresponding
nonlocal problem we consider is written as
\begin{equation}\label{ignition}
-U''-\c U'=\begin{cases} 0 &
\text{ where } U<\theta\\
(U-\theta)(1-\phi *U) &\text{ where } U\geq \theta.
\end{cases}
\end{equation}
Then, one can construct a solution $(c,u)=(c_a,u_a)$ in a bounded
box $[-a,a]$ exactly as in Section \ref{s:tw}, and thus a solution
$(\c,U)$ of \eqref{ignition} as a limit of solutions $(c,u)$. One
can also readily get Proposition \ref{u0iffx0}, which in turn
implies that the solution $u$ solves $-u''-cu'=0$ on $(0,a)$,
$u(0)=\theta$, $u(a)=0$ and therefore becomes explicit on this
interval:
 \begin{equation}\label{explicit}
  u(x)=\frac{-\theta}{e^{ca}-1}+\frac{\theta
  e^{-cx}}{1-e^{-ca}},\;
  \text{ for } 0\leq x \leq a.
  \end{equation}
Assume by contradiction that $\c\leq 0$. Then
  Proposition \ref{propboundaway} $(ii)$, which also directly applies to the ignition case,
implies that there exists $\ep>0$ such that, for $a>0$ large
enough, $u\geq \theta+\ep$
   on $(-\infty,-1/\ep)$, which in turn implies the non triviality
   of $U$. If $\c<0$ then, as $a\to \infty$,
   $$
u'(0)=\frac{-c\theta}{1-e^{-ca}}\to 0,
$$
since $c\to \c<0$. Hence $U'(0)=0$ and then $U\equiv \theta$, a
contradiction. If $\c=0$ then $U$ is a bounded solution of
$-U''=0$ on $(0,\infty)$ such that $U(0)=\theta$, that is $U\equiv
\theta$, a contradiction. As a result $\c>0$. Letting $a\to
\infty$ in \eqref{explicit} yields
$$
U(x)=\theta e^{-\c x},\; \text{ for all } x\geq 0.
$$
To conclude, the behavior \eqref{gauche} as $x\to -\infty$ is
proved as in Proposition \ref{prop:uinfini}.

\bigskip

\noindent \textbf{Acknowledgements.}  M. A. is supported by the
French {\it Agence Nationale de la Recherche} within the project
IDEE (ANR-2010-0112-01).

\end{document}